\documentclass[a4paper]{article}
\usepackage{amsmath}\usepackage{epsf,amsfonts,amsthm}
\usepackage{fontenc,indentfirst, delarray,amsfonts,amsmath,amssymb}

\newcommand{\be}{\begin{equation}}
\newcommand{\ee}{\end{equation}}
\newcommand{\bea}{\begin{eqnarray*}}
\newcommand{\eea}{\end{eqnarray*}}
\newcommand{\w}{\wedge}
\newcommand{\p}{\partial}

\newcommand{\raa}{\rightarrow}

\newcommand{\E}{\ell}
\newcommand{\N}{\mathbb{N}}
\newcommand{\Z}{\mathbb{Z}}
\newcommand{\R}{\mathbb{R}}

\newcommand{\Q}{\mathbb{Q}}

\newcommand{\lp}{\left(}
\newcommand{\rp}{\right)}

\newcommand{\op}[1]{\!\!\mathop{\rm ~#1}\nolimits}

\mathchardef\za="710B  
\mathchardef\zb="710C  
\mathchardef\zg="710D  
\mathchardef\zd="710E  
\mathchardef\zve="710F 
\mathchardef\zz="7110  
\mathchardef\zh="7111  
\mathchardef\zy="7112 
\mathchardef\zi="7113  
\mathchardef\zk="7114  
\mathchardef\zl="7115  
\mathchardef\zm="7116  
\mathchardef\zn="7117  
\mathchardef\zx="7118  
\mathchardef\zp="7119  
\mathchardef\zr="711A  
\mathchardef\zs="711B  
\mathchardef\zt="711C  
\mathchardef\zu="711D  
\mathchardef\zf="711E 
\mathchardef\zq="711F  
\mathchardef\zc="7120  
\mathchardef\zw="7121  
\mathchardef\ze="7122  
\mathchardef\zvy="7123  
\mathchardef\zvw="7124  
\mathchardef\zvr="7125 
\mathchardef\zvs="7126 
\mathchardef\zvf="7127  
\mathchardef\zG="7000  
\mathchardef\zD="7001  
\mathchardef\zY="7002  
\mathchardef\zL="7003  
\mathchardef\zX="7004  
\mathchardef\zP="7005  
\mathchardef\zS="7006  
\mathchardef\zU="7007  
\mathchardef\zF="7008  
\mathchardef\zW="700A  

\begin{document}

\title{Formal Poisson cohomology of \\twisted $r$-matrix induced structures\footnote{The research of N. Poncin was supported by
grant R1F105L10. This author also thanks the Erwin Schr\"odinger
Institute in Vienna for hospitality and support during his visit
in 2006.}}
\author{Mourad Ammar\footnote{University of Luxembourg, Campus
Limpertsberg, Institute of Mathematics, 162A, avenue de la
Fa\"iencerie, L-1511 Luxembourg City, Grand-Duchy of Luxembourg,
E-mail: mourad.ammar@uni.lu}, Norbert Poncin\footnote{University
of Luxembourg, Campus Limpertsberg, Institute of Mathematics,
162A, avenue de la Fa\"iencerie, L-1511 Luxembourg City,
Grand-Duchy of Luxembourg, E-mail:
norbert.poncin@uni.lu}}\maketitle

\newtheorem{rem}{Remark}
\newtheorem{theo}{Theorem}
\newtheorem{prop}{Proposition}
\newtheorem{lem}{Lemma}
\newtheorem{cor}{Corollary}
\newtheorem{ex}{Example}

\begin{abstract}

Quadratic Poisson tensors of the Dufour-Haraki classification read
as a sum of an $r$-matrix induced structure twisted by a (small)
compatible exact quadratic tensor. An appropriate bigrading of the
space of formal Poisson cochains then leads to a vertically
positive double complex. The associated spectral sequence allows
to compute the Poisson-Lichnerowicz cohomology of the considered
tensors. We depict this modus operandi, apply our technique to
concrete examples of twisted Poisson structures, and obtain a
complete description of their cohomology. As richness of Poisson
cohomology entails computation through the whole spectral
sequence, we detail an entire model of this sequence. Finally, the
paper provides practical insight into the operating mode of spectral sequences. \\



\end{abstract}

\section{Introduction}

It is easily seen that any quadratic Poisson tensor of the
Dufour-Haraki classification (DHC), \cite{DH}, reads
\be\zL=\zL_I+\zL_{II}=aY_{23}+bY_{31}+cY_{12}+\zL_{II},\label{MP}\ee
where $a,b,c\in\R,$ where the $Y_i$ are linear, mutually commuting
vector fields ($Y_{ij}=Y_i\w Y_j$), and where $\zL_{II}$ is---as
$\zL_I$---a quadratic Poisson structure. This entails of course
that $\zL_I$ and $\zL_{II}$ are compatible, i.e. that
$[\zL_I,\zL_{II}]=0$, where $[.,.]$ is the Schouten bracket.
Except for structure $10$ of the DHC, where
$\zL_{II}=(3b+1)(y^2-2xz)\p_{23}$
($\p_{23}=\p_{x_2}\p_{x_3}=\p_y\p_z$), the second Poisson
structure is always Koszul-exact, i.e.
$$\zL_{II}=\zP_{\zf}:=(\p_1\zf)\p_{23}+(\p_2\zf)\p_{31}+(\p_3\zf)\p_{12},\quad
\zf\in{\cal S}^3\R^{3*}.$$

In \cite{PX}, P. Xu has proved that any quadratic Poisson tensor
of $\R^{3}$ reads \be\zL=\frac{1}{3}K\w {\cal
E}+\zP_f,\label{PX}\ee where $K$ is the curl of $\zL$, ${\cal E}$
the Euler field, and $f\in{\cal S}^3\R^{3*}$.

In most cases (only cases 9 and 10 of the DHC are exceptional),
term $\zL_I$ of Equation (\ref{MP}), which is twisted by the exact
term $\zL_{II}$ and is---as easily seen---implemented by an
$r$-matrix in the stabilizer $\frak{g}_{\zL}\w\frak{g}_{\zL},$
$\frak{g}_{\zL}=\{A\in\op{gl}(3,\R):[A,\zL]=0\}$, is given by
$$\zL_I=\frac{1}{3}K\w{\cal E}+\zP_{\zl D},$$ where $\zl\in\R^*$ and
$D=\op{det}(Y_1,Y_2,Y_3)$, whereas
$$\zL_{II}=\zP_{\zf}=\zP_{f-\zl D}.$$

Hence, the difference between decompositions (\ref{MP}) and
(\ref{PX}) is that in (\ref{MP}) the biggest possible part of
$\zL$ is incorporated into the $r$-matrix induced structure,
whereas in (\ref{PX}) it is incorporated into the exact structure.

We privilege decomposition (\ref{MP}), since a general computing
technique allows to deal with the cohomology of $\zL_I$,
\cite{MP}, and $\zL_{II}$ vanishes in many cases. In most of the
cases where the small exact tensor $\zL_{II}$ does not vanish, the
decomposition
$$\p_{\zL}:=[\zL,.]=[\zL_I,.]+[\zL_{II},.]=:\p_{\zL_I}+\p_{\zL_{II}}, \quad \p_{\zL_I}^2=\p_{\zL_{II}}^2=\p_{\zL_I}\p_{\zL_{II}}+\p_{\zL_{II}}\p_{\zL_{I}}=0$$
leads to a vertically positive double complex and the
corresponding spectral sequence allows to deduce bit by bit the
cohomology of $\zL$ from that of $\zL_{I}$.\\

In Section $2$, we show how twisted $r$-matrix induced tensors
generate vertically positive double complexes. As richness of
Poisson cohomology entails computation through the whole
associated spectral sequence, we detail a complete model of the
sequence in Section $3$. Section $4$ contains the computation of
the cohomology of tensor $\zL_4$ of the Dufour-Haraki
classification. More precisely, Subsection $4.1$ provides the
second term of the spectral sequence, i.e. the cohomology of the
$r$-matrix induced part $\zL_{4,I}$ of $\zL_4$, which is
accessible to the general cohomological technique developed in
\cite{MP}. After some preliminary work in Subsections $4.2$ and
$4.3$, we are prepared to compute, in Subsection $4.4,$ through
the entire spectral sequence, see Theorem \ref{basicTheo}. As we
aim at the extraction of ``true results'', we are obliged to
detail all the isomorphisms involved in the theory of spectral
sequences and to read our upshots through these isomorphisms.
Hence, in particular, a study of the limiting process in the
sequence and of the reconstruction of the cohomology, precedes, in
Subsection $4.5.1$, the concrete description of the cohomology of
twisted structure $\zL_4$, see Theorem \ref{IndLim} in Subsection
$4.5.2$, and of twisted tensor $\zL_8$, Theorem \ref{IndLim8} in Subsection $5$.\\

The description of the main features of the cohomology of
$r$-matrix induced Poisson structures has been given in \cite{MP}.
The tight relation between Casimir functions and Koszul-exactness
of these Poisson tensors is recalled in Subsection $4.1$, see
Equation (\ref{quasi-exact}) (a generalization can be found in
Subsection $4.3$, see Equation(\ref{extQuasiExact})). Since our
$r$-matrix induced Poisson structures are built with infinitesimal
Poisson automorphisms $Y_i$, see Equation (\ref{MP}), the wedge
products of the $Y_i$ constitute a priori ``privileged'' cocycles.
The associative graded commutative algebra structure of the
Poisson cohomology space now explains part of the cohomology
classes. The second and third term of this cohomology space
contain, in addition to the just mentioned wedge products of
Casimir functions and infinitesimal automorphisms $Y_i$,
non-bounding cocycles the coefficients of which are---in a broad
sense---polynomials on the singular locus of the considered
Poisson tensor. The ``weight in cohomology'' of the singularities
increases with closeness of the Poisson structure to
Koszul-exactness. The appearance of some ``accidental
Casimir-like'' non bounding cocycles completes the
depiction of the main characteristics of the cohomology.\\

If the $r$-matrix induced structure is twisted by an exact
quadratic tensor, the aforementioned spectral sequence constructs
little by little the cohomology of $\zL$ from that of $\zL_I$. In
the examined cases, the basic Casimir $C_I$ of $\zL_I$ is the
first term of the expansion by Newton's binomial theorem of the
basic Casimir $C$ of $\zL$. Beyond the emergence of systematic
conditions on the coefficients of the powers $C^{i},$ $i\in\N,$
and the methodic disappearance of monomials on the singular locus
of $\zL_I$, the main impact on Poisson cohomology of twist
$\zL_{II}$ is the (partial) passage from first term $C_I$ to
complete expansion $C$, a change that takes place gradually for
all powers of these Casimirs, as we compute through the spectral
sequence.

\section{Vertically positive double complex}

\subsection{Definition}

Let $(K,d)$ be a complex, i.e. a differential space, made up by a
graded vector space $K=\oplus_{n\in\N}K^n$ and a differential
$d:K^n\raa K^{n+1}$ that has weight $1$ with respect to this
grading. Assume that each term $K^n$ is itself graded,
$$K^n=\oplus_{r,s\in\N,r+s=n}K^{rs},$$ so that
$K=\oplus_{r,s\in\N}K^{rs}$ is bigraded. We will refer to grading
$K=\oplus_{n\in\N}K^n$ as the diagonal grading. Let
$p,q\in\N,p+q=n$. Differential $d:K^{pq}\raa
\oplus_{r,s\in\N,r+s=n+1}K^{rs}$ induces linear maps
$$d_{ab}:K^{pq}\raa K^{p+a,q+b}\quad (a,b\in\Z,a+b=1),$$ such that $$d=\sum_{a,b\in Z,a+b=1}d_{ab}.$$
If $d_{ab}=0,\forall b<0$ (resp. $d_{ab}=0,\forall a<0$), the
preceding complex is a {\it vertically positive double complex}
(VPDC) (resp. a {\it horizontally positive double complex}
(HPDC)). Vertically positive and horizontally positive double
complexes are {\it semi-positive double complexes}. A complex that
is simultaneously a VPDC and a HPDC is a {\it double complex} (DC)
in the usual sense.

We filter a VPDC (resp. a HPDC) using the {\it horizontal
filtration} (resp. {\it vertical filtration})
$$^h\!K_p=\oplus_{r\in\N,s\ge p}K^{rs}\quad (\mbox{resp.  }^v\!K_p=\oplus_{r\ge p,s\in\N}K^{rs}).$$
These filtrations are compatible (in the usual sense) with the
diagonal grading and differential $d$. Moreover, they are regular,
i.e. $K_p\cap K^n=0,\forall p>n$ (as well for $K_p=^h\!\!\!K_p$ as
for $K_p=^v\!\!\!K_p$), and verify $K_0=K$ and $K_{+\infty}=0$.

The (convergent) spectral sequence (SpecSeq) associated with this
graded filtered differential space is extensively studied below.
Let us stress that in the following we prove several general
results on spectral sequences, which we could not find in
literature. In order to increase the reader-friendliness of our
paper and to avoid scrolling, we chose to give these upshots in
separate subsections that directly precede those where the results
are needed.

\subsection{Application to twisted $r$-matrix induced Poisson structures}

We will now associate a VPDC to twisted $r$-matrix induced Poisson
tensors. Let
$$\zL=\zL_I+\zL_{II}=aY_{23}+bY_{31}+cY_{12}+\zP_{\zf}$$ be as in
Equation (\ref{MP}).

Set $Y_i=\E_{ij}\p_j,$ $\E_{ij}\in\R^{3*}$ (we use the Einstein
summation convention) and $D=\op{det}\E=\op{det}(\E_{ij})\in{\cal
S}^3\R^{3*}$. If $L\in\op{gl}(3,{\cal S}^2\R^{3*})$ is the matrix
of algebraic $(2\times 2)$-minors of $\E$, we have
$\p_i=\frac{L_{ji}}{D}Y_j.$ The formal Poisson cochain space
${\cal P}$ is made up by the $0-$, $1-$, $2-$, and $3-$cochains
\be
C^0=\frac{\zs}{D},C^1=\frac{\zs_1}{D}Y_1+\frac{\zs_2}{D}Y_2+\frac{\zs_3}{D}Y_3,
C^2=\frac{\zs_1}{D}Y_{23}+\frac{\zs_2}{D}Y_{31}+\frac{\zs_3}{D}Y_{12},C^3=\frac{\zs}{D}Y_{123},\label{cochainsYbasis}\ee
where $\zs,\zs_1,\zs_2,\zs_3\in\R[[x_1,x_2,x_3]]$ and where $\zs$,
$\E_{ij}\zs_i$, $L_{ij}\zs_i$ are divisible by $D$ (for any $j$;
$3$-cochains do not generate any divisibility condition). In order
to understand these results, note first that, if ${\cal
L}\in\op{gl}(3,{\cal S}^4\R^{3*})$ denotes the matrix of algebraic
$(2\times 2)$-minors of $L$, we have ${\cal
L}=(\op{det}L)\tilde{L}^{-1}$ and $L=(\op{det}\E)\tilde{\E}^{-1}.$
The last equation entails that $\op{det}L=(\op{det}\E)^2$ and that
$L^{-1}=\frac{1}{\op{det}\E}\tilde{\E}.$ Hence, it follows from
the first equation that ${\cal L}=(\op{det}\E)\E=D\E.$ Let now
$C^2=\zs_1\p_{23}+\zs_2\p_{31}+\zs_3\p_{12}$ be an arbitrary
$2$-cochain. Since its first term reads
\be\zs_1\p_{23}=\frac{\zs_1}{D^2}L_{j2}L_{k3}Y_{jk}=\frac{\zs_1}{D^2}\lp
{\cal L}_{11}Y_{23}+{\cal L}_{21}Y_{31}+{\cal L}_{31}Y_{12}\rp
=\frac{\zs_1}{D}\lp\E_{11}Y_{23}+\E_{21}Y_{31}+\E_{31}Y_{12}\rp,\label{passDelY}\ee
its is clear that any $2$-cochain can be written as announced.
Conversely, the first term of any $2$-vector
$C^2=\frac{\zs_1}{D}Y_{23}+\frac{\zs_2}{D}Y_{31}+\frac{\zs_3}{D}Y_{12}$
reads
$$\frac{\zs_1}{D}Y_{23}=\frac{\zs_1}{D}\E_{2j}\E_{3k}\p_{jk}=\frac{\zs_1}{D}\lp
L_{11}\p_{23}+L_{12}\p_{31}+L_{13}\p_{12}\rp.$$ Thus, such a
$2$-vector $C^2$ is a formal Poisson $2$-cochain if and only if
$L_{ij}\zs_i$ is divisible by D for any $j$. The proofs of the
statements concerning $0$-, $1$-, and $3$-cochains are similar.

Hence, if we substitute the $Y_i$ for the standard basic vector
fields $\p_i$, the cochains assume---roughly speaking---the shape
$\sum f\mathbf{Y}$, where $f$ is a function and $\mathbf{Y}$ is a
wedge product of basic fields $Y_i$. Then the Lichnerowicz-Poisson
coboundary operator $\p_{\zL_I}=[\zL_I,\cdot]$ is just
\be\p_{\zL_I}(f\mathbf{Y})=[\zL_I,f\mathbf{Y}]=[\zL_I,f]\w\mathbf{Y}.\label{cobAdmiss}\ee
More precisely, the coboundary operator associated with $\zL_I$ is
given by \be[\zL_I,C^0]=\nabla C^0,[\zL_I,C^1]=\nabla\w
C^1,[\zL_I,C^2]=\nabla.C^2,\mbox{ and
}[\zL_I,C^3]=0,\label{Lambda1cob}\ee where
$\nabla=\sum_iX_i(\cdot)Y_i$, $X_1=cY_2-bY_3,
X_2=aY_3-cY_1,X_3=bY_1-aY_2,$ and where the {\small RHS} have to
be viewed as notations that give the coefficients of the
coboundaries in the $Y_i$-basis. For instance,
$[\zL_I,C^2]=(\sum_iX_i(\frac{\zs_i}{D}))Y_{123}$.

Of course the formal power series $\zs,\zs_1,\zs_2,\zs_3$ in
Equation (\ref{cochainsYbasis}) read
$$\sum_{J\in\N^3}c_JX^J=\sum_{j_1=0}^{\infty}\sum_{j_2=0}^{\infty}\sum_{j_3=0}^{\infty}c_{j_1j_2j_3}x_1^{j_1}x_2^{j_2}x_3^{j_3}\quad (c_{j_1j_2j_3}\in\R).$$
The degrees $j_1,j_2,j_3\in\N$ and the cochain degree
$c\in\{0,1,2,3\}$ induce a $4$-grading of the formal Poisson
cochain space ${\cal P}$ of polyvector fields with coefficients in
formal power series. Let us emphasize that the degrees $j_i$ are
read in the numerators $\zs$ of the decomposition
$C=\sum\frac{\zs}{D}\mathbf{Y}$. They are tightly related with the
$r$-matrix induced nature of $\zL_I$ and were basic in the method
developed in \cite{MP}. In the following we use the degrees
$r=j_1+j_2+c$ and $s=j_3$ (depending on the considered Poisson
tensor, other degrees could be used, but the preceding ones
encompass the majority of twisted structures) that generate a
bigrading of ${\cal P}$, ${\cal P}=\oplus_{r,s\in\N}{\cal
P}^{rs}.$ When defining the diagonal degree $n=r+s$, we get a
graded space
$${\cal P}=\oplus_{n\in\N}{\cal P}^n,{\cal P}^n=\oplus_{r,s\in\N,r+s=n}{\cal P}^{rs}.$$

We now determine the weights of the coboundary operators
$\p_{\zL_I}$ and $\p_{\zL_{II}}$ with respect to $r$ and $s$.
Actually $D$ is an eigenvector of the basic fields $Y_i$, hence of
the fundamental fields $X_i$, $Y_iD=\zl_iD$, $X_iD=\zm_iD,$
$\zl_i,\zm_i\in\R$. Indeed, since $\zp_{\zl
D}=\zl(\p_1D\,\p_{23}+\p_2D\,\p_{31}+\p_3D\,\p_{12})$, it follows
from Equation (\ref{passDelY}) (take $\zs_j=\zl\p_jD$) (and its
cyclic permutations) that
$$\zp_{\zl D}=\frac{\zl}{D}\lp Y_1D\,Y_{23}+Y_2D\,Y_{31}+Y_3D\,Y_{12}\rp.$$
But $\zp_{\zl D}$ is part of $\zL_I$ and is---more precisely---of
type (\ref{MP}), i.e. reads $$\zp_{\zl D}={\frak l}_1Y_{23}+{\frak
l}_2Y_{31}+{\frak l}_3Y_{12}$$ (${\frak l}_1,{\frak l}_2,{\frak
l}_3\in\R$). Hence, $$Y_iD=\frac{{\frak l}_i}{\zl}D=:\zl_iD,
\forall i\in\{1,2,3\}.$$ In view of Equations
(\ref{cochainsYbasis}) and (\ref{Lambda1cob}), the degrees
$j_1,j_2,j_3$ of the $\zL_I$-coboundary $\p_{\zL_I}C$ of any
cochain $C$ only depend on the values $X_i\lp\frac{\zs}{D}\rp$ of
the fundamental linear fields $X_i$ for an arbitrary formal power
series $\zs=\sum_Jc_JX^J$. Since
$$X_i\lp\frac{\zs}{D}\rp=\sum_Jc_J\frac{1}{D}\lp X_i-\zm_i\op{id}\rp X^J,$$
it is clear that  $\p_{\zL_I}$ preserves the total degree
$\frak{t}=j_1+j_2+j_3$.

In the following, we focus on the first twisted quadratic Poisson
structures that appear in the DHC, i.e. on classes 4, 8, and 11,
see \cite{DH}. Let us recall that
$$\zL_4=ayz\p_{23}+axz\p_{31}+\lp
bxy+z^2\rp\p_{12}=aY_{23}+aY_{31}+bY_{12}+\frac{z^3}{D}Y_{12}=\zL_{4,I}+\zL_{4,II},$$
$$a\neq 0, b\neq 0, Y_1=x\p_1,Y_2=y\p_2,Y_3=z\p_3, D=xyz,$$
$$\zL_8=\lp\frac{a+b}{2}(x^2+y^2)\pm z^2\rp\p_{12}+axz\p_{23}+ayz\p_{31}=aY_{23}+\frac{a+b}{2}Y_{12}\pm\frac{z^3}{D}Y_{12}=\zL_{8,I}+\zL_{8,II},$$
$$a\neq 0, b\neq 0, Y_1=x_1\p_1+x_2\p_2,Y_2=x_1\p_2-x_2\p_1,Y_3=x_3\p_3, D=(x^2+y^2)z,$$
$$\zL_{11}=\lp ax^2+bz^2\rp\p_{12}+(2a+1)xz\p_{23}=Y_{23}+aY_{12}+b\frac{z^3}{D}\lp (3a+1)Y_{12}+Y_{23}\rp=\zL_{11,I}+\zL_{11,II},$$
$$a\neq\frac{-1}{3},b\neq 0, Y_1={\cal E}, Y_2=x\p_2, Y_3=(3a+1)z\p_3, D=(3a+1)x^2z.$$
Owing to the above remarks, it is obvious that $\p_{\zL_{i,I}},
i\in\{4,8,11\},$ preserves the partial degree $\frak{p}=j_1+j_2$
(and, as aforementioned, the total degree $\frak{t}$). Hence, its
weight with respect to $(r,s)$ is $(1,0)$:
$$d':=d_{10}:=\p_{\zL_{i,I}}:{\cal P}^{rs}\raa{\cal P}^{r+1,s}\quad (i\in\{4,8,11\})$$
(dependence on $i$ omitted in $d'$ and $d_{10}$).

As for the weight of $\p_{\zL_{i,II}}, i\in\{4,8,11\},$ with
respect to $(r,s)$, let us first recall that, if $f$ and $g$ are
some functions, and if $\mathbf{X}$ and $\mathbf{Y}$ denote wedge
products of $Y_1,Y_2,Y_3$ with (non-shifted) degrees $\za$ and
$\zb$ respectively, we have
\be[f\mathbf{X},g\mathbf{Y}]=f[\mathbf{X},g]\w\mathbf{Y}+(-1)^{\za\zb-\za-\zb}g[\mathbf{Y},f]\w\mathbf{X}.\label{bracketFuncWedge}\ee
Of course, the {\small RHS} of the preceding equation is a linear
combination of terms of the type $fY_i(g)\mathbf{Z}$ or
$gY_i(f)\mathbf{Z}$, where $\mathbf{Z}$ is a wedge product of
$Y_1,Y_2,Y_3$ of degree $\za+\zb-1$. It follows that
$\p_{\zL_{i,II}}C^c$, $i\in\{4,8,11\}$, $C^c\in{\cal P}$, is a
formal series of terms of the type
$$[\frac{z^3}{D}\mathbf{X},\frac{X^J}{D}\mathbf{Y}].$$ Any such term is a
linear combination of terms of the type
$$\frac{z^3}{D}Y_i\lp\frac{X^J}{D}\rp\mathbf{Z}\quad\mbox{and}\quad\frac{X^J}{D}Y_i\lp\frac{z^3}{D}\rp\mathbf{Z}.$$
As $D$ is an eigenvector of $Y_i$, this entails that coboundary
$\p_{\zL_{i,II}}C^c$ has the form
$$\p_{\zL_{i,II}}C^c=\sum\frac{\sum_K c_K X^K}{D^2}\mathbf{Z},$$
where in each term $k_1+k_2=j_1+j_2$ and $k_3=j_3+3,$ and where
the degree of wedge product $\mathbf{Z}$ is $\za+\zb-1=c+1$. When
dividing the preceding numerators by $D$ (see above), we find that
the weight of $\p_{\zL_{i,II}}$ with respect to $(r,s)$ is
$(-1,2):$
$$d'':=d_{-12}:=\p_{\zL_{i,II}}:{\cal P}^{rs}\raa{\cal P}^{r-1,s+2}\quad (i\in\{4,8,11\})$$ (dependence on $i$ omitted in
$d''$ and $d_{-12}$).

Finally, $({\cal P},\p_{\zL_i})$, $i\in\{4,8,11\}$, endowed with
the previously mentioned gradings $${\cal P}=\oplus_{n\in\N}{\cal
P}^n,{\cal P}^n=\oplus_{r,s\in\N,r+s=n}{\cal P}^{rs}$$ and the
differential
$$d:=\p_{\zL_i}=\p_{\zL_{i,I}}+\p_{\zL_{i,II}}=d'+d''=d_{10}+d_{-12},$$ is a
VPDC. We will compute the cohomology $H(\zL_i)=H({\cal P},d)$
using the SpecSeq associated with this VPDC (see above).

\section{Model of the spectral sequence associated with a
VPDC}\label{ModSpecSeq}

As mentioned above, a VPDC, a HPDC, and a DC can canonically be
viewed as regular filtered graded differential spaces. Hence, a
SpecSeq (two, for any DC) is associated with each one of these
complexes.

In order to introduce notations, let us recall that, if $(K,d,
K_p,K^n)$ is any (regular, i.e. $K_p\cap K^n=0,\forall p>n$)
filtered (subscripts) graded (superscripts) differential space (in
our work $p$ and $n$ can be regarded as positive integers), the
associated SpecSeq $(E_r,d_r)$ ($r\in\N$) is defined by
$$E^{pq}_r=Z_r^{pq}/(Z_{r-1}^{p+1,q-1}+B_{r-1}^{pq}),$$
where $Z_r^{pq}=K_p\cap d^{-1}K_{p+r}\cap K^{p+q}$ and
$B_r^{pq}=K_p\cap dK_{p-r}\cap K^{p+q}$ are the spaces of ``weak
cocycles'' and ``strong coboundaries'' of order $r$ in $K_p\cap
K^{p+q}$, and $$d_r:E_r^{pq}\ni[\frak{z}_r^{pq}]_{E_r^{pq}}\raa
[d\frak{z}_r^{pq}]_{E_r^{p+r,q+1-r}}\in E_r^{p+r,q+1-r}.$$ In the
following, we also use the vector space isomorphism
$$\zs_r:E_{r+1}^{pq}\raa H^{pq}(E_r,d_r),$$ which assigns to each
$[\frak{z}_{r+1}^{pq}]_{E_{r+1}^{pq}}$, $\frak{z}_{r+1}^{pq}\in
Z_{r+1}^{pq}\subset Z_r^{pq},$ the $d_r$-cohomology class
$[[\frak{z}_{r+1}^{pq}]_{E_r^{pq}}]_{d_r}$,
$[\frak{z}_{r+1}^{pq}]_{E_r^{pq}}\in E_r^{pq}\cap\op{ker}d_r$. 
For more detailed results on spectral sequences, we refer the
reader to \cite{JMC}, \cite{RG}, \cite{CE}, \cite{IV1}, ... In
these monographs, a model for the SpecSeq associated with a (HP)DC
is partially depicted up to $r=2$. It is well-known that spectral
sequences are particularly easy to use, if many spaces $E_2^{pq}$
(or $E_r^{pq}$ $(r>2)$) vanish. Due to richness of Poisson
cohomology, this lacunary phenomenon is less pronounced in our
setting. Since we have thus to compute through the whole SpecSeq,
we need the complete description of the entire model of the
SpecSeq $(E_r,d_r)$ ($r\in\N$) associated with a VPDC.\\

So consider an arbitrary VPDC and let $G^{pq}(K)$ ($p,q\in\N$) be
the term of degree $(p,q)$ of the bigraded space associated with
the filtered graded space $K$. It is clear that the mapping
$$I_0: E_0^{pq}=K_p\cap K^{p+q}/K_{p+1}\cap K^{p+q}=G^{pq}(K)\ni
[\frak{z}_0^{pq}=\sum_{i=0}^q z^{q-i,p+i}]_{E_0^{pq}}\raa
z^{qp}\in K^{qp},$$ where $z^{rs}$ (as well as---in the
following---all Latin characters with double superscript) is an
element of $K^{rs}$ (whereas German Fraktur characters with double
superscript, such as $\frak{z}_0^{pq}$, do not refer to the
bigrading of $K$), is an isomorphism of bigraded vector spaces
(i.e. a vector space isomorphism that respects the bigrading). It
is easily seen that, when reading $d_0$ through this isomorphism,
we get the compound map
$$\overline{d}_0=I_0d_0I_0^{-1}=d_{10}.$$

Thus $I_0: (E_0,d_0)\raa (K,\overline{d}_0)$ is an isomorphism
between bigraded differential spaces, and induces an isomorphism
$$I_{0\sharp}:H^{pq}(E_0,d_0)\ni [[\frak{z}_0^{pq}=\sum_{i=0}^q z^{q-i,p+i}]_{E_0^{pq}}]_{d_0}\raa
[z^{qp}]_{\overline{d}_0}\in
H^{pq}(K,\overline{d}_0)=:\,^0\!H^{pq}(K)=\,^0\!H^q(K^{*p})$$ of
bigraded vector spaces, where the last space is the $q$-term of
the cohomology space of $(K^{*p},\overline{d}_0=d_{10})$. Hence
the bigraded vector space isomorphism
$$I_1=I_{0\sharp}\zs_0:E_1^{pq}\ni [\frak{z}_1^{pq}=\sum_{i=0}^q z^{q-i,p+i}]_{E_1^{pq}}\raa [[\frak{z}_1^{pq}]_{E_0^{pq}}]_{d_0}\raa
[z^{qp}]_{\overline{d}_0}\in\,^0\!H^q(K^{*p}).$$ We now again
verify straightforwardly that differential $d_1$ read on model
$^0\!H(K)$ is induced by $d_{01}$, i.e. that
$$\overline{d}_1=I_1d_1I_1^{-1}=d_{01\sharp}.$$

Finally, $$I_2=I_{1\sharp}\zs_1:E_2^{pq}\ni
[\frak{z}_2^{pq}=\sum_{i=0}^q z^{q-i,p+i}]_{E_2^{pq}}\raa
[[\frak{z}_2^{pq}]_{E_1^{pq}}]_{d_1}\raa
[[z^{qp}]_{\overline{d}_0}]_{\overline{d}_1}\in
^1\!\!\!H^p(^0\!H^q(K))$$ is an isomorphism of bigraded vector
spaces. As for the sense of the last space, note that
$(^0\!H^q(K)=\oplus_p\,^0\!H^q(K^{*p}),\overline{d}_1)$ is a
complex. Observe now that the inverse $I_2^{-1}$ is less
straightforward than $I_0^{-1}$ and $I_1^{-1}$. Indeed, if
$[[z^{qp}]_{\overline{d}_0}]_{\overline{d}_1}\in
^1\!\!\!H^p(^0\!H^q(K))$, representative $z^{qp}$ is generally not
a member of $Z_2^{pq}$. However, since the considered class makes
sense,
$$\begin{array}{l}d_{10}z^{qp}=0\\d_{01}z^{qp}+d_{10}z^{q-1,p+1}=0,\end{array}$$
where $z^{q-1,p+1}\in K^{q-1,p+1}$. Thus,
$\frak{z}_2^{pq}:=z^{qp}+z^{q-1,p+1}\in Z_2^{pq}$ and
$$I_2^{-1}[[z^{qp}]_{\overline{d}_0}]_{\overline{d}_1}=[\frak{z}_2^{pq}]_{E_2^{pq}}.$$
So
$$\overline{d}_2[[z^{qp}]_{\overline{d}_0}]_{\overline{d}_1}=I_2[d\frak{z}_2^{pq}]_{E_2^{p+2,q-1}}=[[d_{-12}z^{qp}+d_{01}z^{q-1,p+1}]_{\overline{d}_0}]_{\overline{d}_1}.$$
The preceding results extend those given in \cite{IV1} (for a
HPDC). They can easily be adapted to the most frequently
encountered situations where only some terms
$d_{ab}$ of $d$ do not vanish.\\

In the following, we complete the description of the SpecSeq
associated with a VPDC, assuming that $d=d_{10}+d_{-12}:=d'+d''$.
This hypothesis entails that $d'^2=d''^2=d'd''+d''d'=0$, i.e. that
$d'$ and $d''$ are two anticommuting differentials. Hereafter, we
denote by $^r\!H(.)$ ($r\in\N$) the cohomology of differential
$\overline{d}_{2r}$ and by $[.]_r$ the corresponding classes.
Moreover, we will deal with strongly triangular systems of type
$$\begin{array}{ll}d'z^{qp}=0&(\frak{E}_0)\\d''z^{qp}+d'z^{q-2,p+2}=0&(\frak{E}_1)\\\ldots\\d''z^{q-2(k-2),p+2(k-2)}+d'z^{q-2(k-1),p+2(k-1)}=0.&(\frak{E}_{k-1})\end{array}$$
Note that when solving such a system, we prove at each stage that
some $d'$-cocycle is actually a $d'$-coboundary. We refer to this
kind of system using the notation $S(z^{qp};k)$ or
$S(k;z^{q-2(k-1),p+2(k-1)})$ depending on the necessity to
emphasize the first or the last unknown or entry of an ordered
solution.

\begin{prop}\label{modelSpecSeq} The spectral sequence associated to a VPDC with differential
$d=d_{10}+d_{-12}=d'+d''$ admits the following model. The model of
$E_0$, isomorphisms $I_0$ and $I_0^{-1}$, and
differential $\overline{d}_0$ are the same as above. For any $r\in\{1,2,\ldots\}$,\\

(i) The map
$$I_{2r-1}:E_{2r-1}^{pq}\ni[\frak{z}_{2r-1}^{pq}=\sum_{i=0}^qz^{q-i,p+i}]_{E_{2r-1}^{pq}}\raa
[[[z^{qp}]_0]_1\ldots]_{r-1}\in\, ^{r-1}\!H^{pq}(^{r-2}\!H(\ldots
(^0\!H(K))))$$ is a bigraded vector space isomorphism. Its inverse
$I_{2r-1}^{-1}$ associates to any {\small RHS}-class the {\small
LHS}-class with representative
$\frak{z}_{2r-1}^{pq}=\sum_{i=0}^{r-1}z^{q-2i,p+2i}$, where
$(z^{qp},\ldots,z^{q-2(r-1),p+2(r-1)})$ is any solution of system
$S(z^{qp};r)$. Furthermore, $\overline{d}_{2r-1}=0.$\\

(ii) The model of $E_{2r}^{pq}$ and the corresponding isomorphisms
$I_{2r}$ and $I_{2r}^{-1}$ coincide with those pertaining to
$E_{2r-1}^{pq}$. Moreover,
\be\overline{d}_{2r}[[[z^{qp}]_0]_1\ldots]_{r-1}=[[[d''z^{q-2(r-1),p+2(r-1)}]_0]_1\ldots]_{r-1},\label{cob2r}\ee
where $z^{q-2(r-1),p+2(r-1)}$ is the last entry of an arbitrary
solution of $S(z^{qp};r)$.
\end{prop}

{\it Proof}. It is easier to prove an extended version of
Proposition \ref{modelSpecSeq}. Indeed, let us complete
assertions (i) and (ii) by item\\

(iii) Existence (resp. vanishing) of a class
$[[[z^{qp}]_0]_1\ldots]_{r-1}$ is equivalent with existence of at
least one solution of system $S(z^{qp};r)$ (resp. with existence
of $z^{q-1,p}$ and of $z_i^{q+1,p-2}$, $i\in\{1,\ldots,r-1\}$,
which induce systems $S(i;z_i^{q+1,p-2})$ with solution, such that
$$z^{qp}+d'z^{q-1,p}+d''\sum_{i=1}^{r-1}z_i^{q+1,p-2}=0.)$$

The proof is by induction on $r$. Observe first that the
assertions are valid for $r=1$ (see above). Assume now that all
items hold for $r\in\{1,\ldots,\E-1\}$. Proceeding as above, we
easily show that $I_{2\E-1}:=I_{2(\E-1)\sharp}\zs_{2(\E-1)}$ is
the appropriate bigraded vector space isomorphism. In order to
determine $I_{2\E-1}^{-1}$, take any {\small RHS}-class
$[[[z^{qp}]_0]_1\ldots]_{\E-1}$.\\

Let us first prove assertion (iii). Existence of class
$[[[z^{qp}]_0]_1\ldots]_{\E-1}$ is equivalent with existence of
class $[[[z^{qp}]_0]_1\ldots]_{\E-2}$ (itself equivalent to
existence of at least one solution
$$z^{q-2j,p+2j}\quad (0\le j\le\E-2)$$ for $S(z^{qp};\E-1)$, by
induction) and condition
$$\overline{d}_{2(\E-1)}[[[z^{qp}]_0]_1\ldots]_{\E-2}=0.$$
Using the induction assumptions, we see that the last condition is
equivalent, first with
$$[[[d''z^{q-2(\E-2),p+2(\E-2)}]_0]_1\ldots]_{\E-2}=0,$$ then with
existence of $$z^{q-2(\E-1),p+2(\E-1)}$$ and
$z_i^{q-2(\E-2),p+2(\E-2)}$ ($1\le i\le\E-2$), which implement
systems $S(i;z_i^{q-2(\E-2),p+2(\E-2)})$ with solution, say
$$z_i^{q-2j,p+2j}\quad (1\le \E-i-1\le j\le\E-2),$$ such that
\be d''\lp
z^{q-2(\E-2),p+2(\E-2)}+\sum_{i=1}^{\E-2}z_i^{q-2(\E-2),p+2(\E-2)}\rp+d'z^{q-2(\E-1),p+2(\E-1)}=0.\label{lastEqS}\ee
Assume now that all this holds and define new $z^{q-2j,p+2j}$
($0\le j\le \E-1$). For each $j$, take just the sum of the old
$z^{q-2j,p+2j}$ and of all existing $z_i^{q-2j,p+2j}$. These new
$z^{q-2j,p+2j}$ form a solution of $S(z^{qp};\E)$. Note first that
for $j\in\{0,\E-1\}$, the old and new $z^{q-2j,p+2j}$ coincide.
Hence, the last equation $(\frak{E}_{\E-1})$ of $S(z^{qp};\E)$ is
nothing but Equation (\ref{lastEqS}). Moreover, it is easily
checked that Equations $(\frak{E}_{\E-2}),\ldots,(\frak{E}_0)$ are
also verified. Conversely, if $S(z^{qp};\E)$ has a solution, the
successive classes
$[z^{qp}]_0,[[z^{qp}]_0]_1,\ldots,[[[z^{qp}]_0]_1\ldots]_{\E-1}$
are actually defined. It suffices to note that
$\overline{d}_0z^{qp}=0$ and that, by induction,
$$\overline{d}_{2r}[[[z^{qp}]_0]_1\ldots]_{r-1}=[[[d''z^{q-2(r-1),p+2(r-1)}]_0]_1\ldots]_{r-1}=-[[[d'z^{q-2r,p+2r}]_0]_1\ldots]_{r-1}=0,$$
for any $r\in\{1,\ldots,\E-1\}$.

As for the second part of (iii), note that a class
$[[[z^{qp}]_0]_1\ldots]_{\E-1}$ vanishes if and only if there is
$z_{\E-1}^{q+2(\E-1)-1,p-2(\E-1)}$ that generates a system
$S(z_{\E-1}^{q+2(\E-1)-1,p-2(\E-1)};\E-1)$ with solution, say
$$z_{\E-1}^{q+2(\E-j-1)-1,p-2(\E-j-1)}\quad (0\le j\le \E-2),$$ such that
$$[[[z^{qp}]_0]_1\ldots]_{\E-2}=-\overline{d}_{2(\E-1)}[[[z_{\E-1}^{q+2(\E-1)-1,p-2(\E-1)}]_0]_1\ldots]_{\E-2}
=-[[[d''z_{\E-1}^{q+1,p-2}]_0]_1\ldots]_{\E-2}.$$ But, by
induction, $[[[z^{qp}+d''z_{\E-1}^{q+1,p-2}]_0]_1\ldots]_{\E-2}=0$
if and only if there are $z^{q-1,p}$ and $z_i^{q+1,p-2}$ ($1\le
i\le \E-2$), which induce systems $S(i;z_i^{q+1,p-2})$ with
solution, such that
$$z^{qp}+d''z_{\E-1}^{q+1,p-2}+d''\sum_{i=1}^{\E-2}z_i^{q+1,p-2}+d'z^{q-1,p}=0.$$
Hence the conclusion.\\

We now revert to items (i) and (ii). For any {\small RHS}-class
$[[[z^{qp}]_0]_1\ldots]_{\E-1}\in\;
^{\E-1}\!H^{pq}(^{\E-2}\!H(\ldots (^0\!H(K))))$, the corresponding
system $S(z^{qp};\E)$ admits, as just explained, at least one
solution $z^{q-2j,p+2j}$ ($0\le j\le \E-1$). Set
$$\frak{z}_{2\E-1}^{pq}:=\sum_{j=0}^{\E-1}z^{q-2j,p+2j}.$$ As
$d\frak{z}_{2\E-1}^{pq}=d''z^{q-2(\E-1),p+2(\E-1)}\in
K^{q-2\E+1,p+2\E}$, we see that $\frak{z}_{2\E-1}^{pq}\in
Z_{2\E-1}^{pq}= K_p\cap d^{-1}K_{p+2\E-1}\cap K^{p+q}$. Hence
$I_{2\E-1}^{-1}$. As
$d_{2\E-1}[\frak{z}_{2\E-1}^{pq}]_{E_{2\E-1}^{pq}}\in
E_{2\E-1}^{p+2\E-1,q-2\E+2}$, it is clear that
$\overline{d}_{2\E-1}=0.$ Thus, the statement concerning the model
of $E_{2\E}^{pq}$ and the isomorphisms $I_{2\E}$ and
$I_{2\E}^{-1}$ is obvious. Finally, as
$d_{2\E}[\frak{z}_{2\E}^{pq}]_{E_{2\E}^{pq}}\in
E_{2\E}^{p+2\E,q-2\E+1},$ we get
$$\overline{d}_{2\E}[[[z^{qp}]_0]_1\ldots]_{\E-1}=[[[d''z^{q-2(\E-1),p+2(\E-1)}]_0]_1\ldots]_{\E-1}.\quad\rule{1.5mm}{2.5mm}$$

{\bf Remark}. Result (\ref{cob2r}) can be rephrased as
$\overline{d}_{2r}=\lp(-1)^{r-1}d''(d'^{-1}d'')^{r-1}\rp_{\sharp}$,
for any $r\in\{1,2,\ldots\}$.

\section{Formal cohomology of Poisson tensor $\zL_4$}

As aforementioned, we use the just depicted SpecSeq associated
with the above detailed VPDC implemented by the twisted $r$-matrix
induced Poisson structure $\zL_4$.

\subsection{Computation of the second term of the SpecSeq}\label{structure}

In this section, we give the second term $E_2\simeq\, ^0\!H(\cal
P)$ of the SpecSeq. Note that $^0\!H(\cal P)$ is the formal
Poisson cohomology of $\overline{d}_0=d'=d_{10}=\p_{\zL_{4,I}}$.
As already elucidated in the Introduction, we came up with
decomposition (\ref{MP}), since the cohomology of $\p_{\zL_I}$ is
always accessible by the technique proposed in \cite{MP}. Hence,
cohomology space $^0\!H(\cal P)$ can be obtained (quite
straightforwardly) by this modus operandi. Let us emphasize that
our results are in accordance, as well with similar upshots in
\cite{PM3}, as with our comments in \cite{MP}, regarding the tight
relation between Casimir functions and Koszul-exactness or
``quasi-exactness'', the appearance of ``accidental Casimir-like''
non bounding cocycles, and the increase of the ``weight in
cohomology'' of the singularities, with closeness of the
considered Poisson structure to Koszul-exactness.

If $\frac{b}{a}\in\Q^*,$ we denote by $(\zb,\za)\sim(b,a)$,
$\za\in\N^*$, the irreducible representative of $\frac{b}{a}$.
Remember that, see \cite{MP}, for $\frac{b}{a}\in\Q^*_+,$ a
quasi-exact structure
\be\zL=a\p_1(pq)\p_{23}+a\p_2(pq)\p_{31}+b\p_3(pq)\p_{12},\label{quasi-exact}\ee
$p=p(x,y),$ $q=q(z),$ exhibits the basic Casimir $p^{\za}q^{\zb}$.
Furthermore, we set $D=xyz,$ $D'=xy$, and write ${\cal
A}_{\za}Y_3$, $\za\in\N^*$, instead of
$D'^{\za}z^{-1}Y_3=D'^{\za}\p_3$, and $\oplus_{ij}\ldots Y_{ij}$
instead of $\ldots Y_{23}+\ldots Y_{31}+\ldots Y_{12}$. Remark
also that the algebra of polynomials of the algebraic variety of
singularities of $\zL_{4,I}$ is
$\R[[x]]\oplus\R[[y]]\oplus\R[[z]]$, where it is understood that
term $\R$ is considered only once.

The following proposition is now almost obvious.
\begin{prop}\label{cohoZL4I}$\quad$

\begin{enumerate}
\item If $\frac{b}{a}\in\Q^*_+,$ the algebra of
$\zL_{4,I}$-Casimirs is $\op{Cas}(\zL_{4,I})=\oplus_{i\in\N}\R
D'^{\za i}z^{\zb i}$ and the cohomology space $^0\!H({\cal P})$ is
given by \begin{eqnarray*} E_2\simeq\,^0\!H({\cal P})
&=&\op{Cas}(\Lambda _{4,I})\oplus\bigoplus_{i}\op{Cas}(\Lambda
_{4,I})Y_{i}\oplus \bigoplus_{ij}\op{Cas}(\Lambda
_{4,I})Y_{ij}\oplus
\op{Cas}(\Lambda _{4,I})Y_{123}\\
&&\oplus \R[[z]]\partial _{12}\oplus\R[[z]]\partial _{123}\oplus
\left\{
\begin{array}{l}
\R[[x]]\p_{23}\oplus\R[[y]]\p_{31}\oplus
(\R[[x]]\oplus\R[[y]])\partial _{123},\text{if }b=a\\
0,\text{otherwise}
\end{array}
\right.
\end{eqnarray*}

\item If $\frac{b}{a}\in\R^*\setminus\Q^*_+,$ we have
$\op{Cas}(\Lambda _{4,I})=\R$ and
\begin{eqnarray*}
E_2\simeq\,^0\!H({\cal P}) &=&\op{Cas}(\Lambda
_{4,I})\oplus\bigoplus_{i}\op{Cas}(\Lambda _{4,I})Y_{i}\oplus
\bigoplus_{ij}\op{Cas}(\Lambda _{4,I})Y_{ij}\oplus
\op{Cas}(\Lambda _{4,I})Y_{123}\\&&\oplus \left\{
\begin{array}{l}
\R{\cal A}_{\za}Y_{3}\oplus {\cal A}_{\za}(\R Y_{23}+\R
Y_{31})\oplus \R{\cal A}_{\za}Y_{123},\text{ if
}(-1,\za)\sim(b,a)\\0,\text{otherwise}
\end{array}
\right.\\&&\oplus \R[[z]]\partial _{12}\oplus\R[[z]]\partial
_{123}
\end{eqnarray*}
\end{enumerate}
\end{prop}

{\bf Remark}. Due to the properties---used below---of the
preceding (non bounding) $\zL_{4,I}$-cocycles, we classify these
representatives as follows:

\begin{enumerate}
\item Representatives of type 1: All cocycles with cochain degree
$0$, the $1-$ and $2-$cocycles that contain a Casimir (maybe the
accidental Casimir ${\cal A}_{\za}$), except cocycles
$\op{Cas}(\zL_{4,I})Y_{12}$

\item Representatives of type 2: All $3-$cocycles, all cocycles
with singularities, and cocycles $\op{Cas}(\Lambda _{4,I})Y_{12}$
\end{enumerate}

\subsection{Prolongable systems $\mathbf{S(z^{qp};r)}$}

Since computation through the whole SpecSeq will shape up as
inescapable, we need the below corollary of Proposition
(\ref{modelSpecSeq}). It allows to short-circuit the process of
computing the successive terms of the sequence. Let us specify
that in the following a system of representatives of a space of
classes is made up by representatives that are in $1$-to-$1$
correspondence with the considered classes.

\begin{cor}\label{Corollary}
If, for some fixed $r\in\N^*$, all the classes
$[[[z^{qp}]_0]_1\ldots]_{r-1}$ in model space
$^{r-1}\!H(^{r-2}\!H(\ldots \linebreak^0\!H(K))),$ appendant on a
SpecSeq associated to a VPDC with differential
$d=d_{10}+d_{-12}=d'+d''$, give rise to an enlarged system
$S(z^{qp};s)$ with solution, for some fixed $s\ge r$, the
following upshots hold:

\begin{enumerate}
\item All the differentials $\overline{d}_{2r-1},
\overline{d}_{2r}, \ldots,\overline{d}_{2s-1}$ vanish \item
Differential $\overline{d}_{2s}$ is defined by
$\overline{d}_{2s}[[[z^{qp}]_0]_1\ldots]_{r-1}=[[[d''z^{q-2(s-1),p+2(s-1)}]_0]_1\ldots]_{r-1}$
\item Any system $(z^{qp})$ of representatives of
$^{r-1}\!H(^{r-2}\!H(\ldots ^0\!H(K)))$ is in $1$-to-$1$
correspondence with the system
$(\frak{z}_{2s}^{pq}:=\sum_{k=0}^{s-1}z^{q-2k,p+2k})$ of
representatives of $E_{2s}$
\end{enumerate}
\end{cor}

{\it Proof}. Induction on $s$. \rule{1.5mm}{2.5mm}

\subsection{Forecast}

In order to increase readability of our paper, some intuitive
advisements are necessary.

The basic idea of the theory of spectral sequences is that
computation of the successive terms $E_{r}\simeq
H(E_{r-1},d_{r-1})$ ($r\in\N^*$) allows to detect their inductive
limit $E_{\infty}$, which---for a convergent sequence---is
isomorphic with the graded space $G(H)$ associated to the
sought-after filtered cohomology space $H$. We then hope to be
able to reconstruct this filtered space $H$ from the corresponding
graded space $G(H)$. Let us recall that space $H$ is of course the
cohomology of the filtered graded differential space associated
with the SpecSeq. Hence, in our case, $H=H(\zL_4)$. It is clear
that the successive cohomology computations take place on the
concrete model side. To determine $H$, we have to pull our results
back to the theoretical side, and more precisely to read them
through the numerous isomorphisms involved.

Actually the application of spectral sequences presented in this
work, provides a beautiful insight into the operating mode of a
SpecSeq. Since---roughly spoken---the ``weak cocycle condition''
in the definition of $Z_r^{pq}$ (resp. the ``strong coboundary
condition'' in the definition of $B_r^{pq}$) converges to the
usual cocycle condition (resp. the usual coboundary condition), we
understand that, when passing from one estimate $E_{r-1}$ of $H$
to the next approximation $E_r$, we obtain an increasing number of
conditions on our initial weak non bounding cocycles of $E_2$ and
an increasing number of bounding cocycles. Moreover, when we
compute through the SpecSeq, the aforementioned pullbacks, see
Proposition (\ref{modelSpecSeq}), add up solutions of crescive
systems,
$$\frak{z}_{2r}^{pq}=z^{qp}+\sum_{k=1}^{r-1}z^{q-2k,p+2k}.$$

The next remarks aim at anticipation of these systems. The reader
is already familiar with Casimirs of exact and quasi-exact
structures. When taking an interest in slightly more general
quasi-exact tensors,
\be\zL=a\p_1((p+r)q)\p_{23}+a\p_2((p+r)q)\p_{31}+b\p_3((p+r)q)\p_{12},\label{extQuasiExact}\ee
$a,b\in\R^*$, $p=p(x,y),$ $q=q(z),$ $r=r(z),$ it is natural to ask
which polynomials of the type $(p+c\,r)^nq^m,$ $c\in\R,$
$n,m\in\N,$ $(n,m)\neq (0,0),$ are Casimir functions. It is easily
checked that structure $\zL_4$ has this form and that the Casimir
conditions read $am=bn$ and $3bn=ca(2n+m).$ So, for
$\frac{b}{a}\in\Q^*_+,$ the basic Casimir $C$ of $\zL_4$ and its
powers $C^{i}$, $i\in\N$, are given by
$$C^{i}=(p+\frac{3b}{2a+b}r)^{\za\,i}q^{\zb\,i}=(D'+\frac{z^2}{2a+b})^{\za\,i}z^{\zb\,i}
=D'^{\za\,i}z^{\zb\,i}+\sum_{k=1}^{\za\,i}\frac{\complement^k_{\za\,i}}{(2a+b)^k}D'^{\za\,i-k}z^{\zb\,i+2k}.$$
These powers $C^{i}$ (non bounding cocycles of $H=H(\zL_4)$) will
be obtained---while we compute through the SpecSeq---from those,
$D'^{\za\,i}z^{\zb\,i}$, of the Casimir of $\zL_{4,I}$ (non
bounding cocycles of $E_2\simeq\,^0\!H({\cal P})$). Hence, the
above-quoted solutions and corresponding systems
$S(D'^{\za\,i}z^{\zb\,i},\za\,i+1)$.

\subsection{Computation through the SpecSeq}

In view of the preceding awareness, it is natural to set
$$Z^{q_{ic}-2k,p_i+2k}=\frac{\complement^k_{\za i}}{(2a+b)^k}D'^{\za i-k}z^{\zb i+2k}\begin{cases}A_{ik}\\B_{ik}Y_1+C_{ik}Y_2+D_{ik}Y_3\quad,\\E_{ik}Y_{23}+F_{ik}Y_{31}\end{cases}$$
where $k\in\{0,1,\ldots,\za\, i\}$ and
$A_{ik},B_{ik},C_{ik},D_{ik},E_{ik},F_{ik}\in\R$. More precisely,
if $\frac{b}{a}\in\Q^*_+,$ we have $(b,a)\sim(\zb,\za),$
$\za,\zb\in\N^*,$ and we ask that $i\in\N$, if
$\frac{b}{a}\in\R^*\setminus\Q^*_+$, we choose $i=0$, and if
moreover $(b,a)\sim(\zb,\za)=(-1,\za)$, $\za\in\N^*$, we also
accept the value $i=1$, but add the conditions
$A_{10}=B_{10}=C_{10}=0$. We define
$(q_{ic},p_i):=(2\za\,i+2+c,\zb\,i+1),$ where $c\in\{0,1,2\}$
denotes the cochain degree, so that the double superscript in the
{\small LHS} is the bidegree $(r,s)=(j_1+j_2+c,j_3)$ of the
{\small RHS}.

Observe that the $Z^{q_{ic},p_i}$ are exactly the representatives
of type 1 of the classes of $E_2\simeq\,^0\!H({\cal P}).$

\begin{lem}\label{Lemma}
For any admissible exponent $i$ and any cochain degree
$c\in\{0,1,2\}$, the cochains $Z^{q_{ic}-2k,p_i+2k}$,
$k\in\{0,1,\ldots,\za\, i\},$ constitute a solution of system
$S(Z^{q_{ic},p_i};\za\,i+1)$, if and only if, for any
$k\in\{0,1,\ldots,\za\, i-1\},$
$$A_{i,k+1}=A_{ik},\;\mbox{if}\;\;c=0,\quad (C_0)$$ $$B_{i,k+1}+C_{i,k+1}=\frac{\lp\za i-k+1\rp \lp B_{ik}+C_{ik}\rp-2D_{ik}}{\za i-k}\quad\mbox{and}
\quad D_{i,k+1}=D_{ik},\;\mbox{if}\;\;c=1,\quad (C_1)$$
$$E_{i,k+1}-F_{i,k+1}=\frac{\za i-k+1}{\za i-k}\lp E_{ik}-F_{ik}\rp,\;\mbox{if}\;\;c=2.\quad (C_2)$$ Furthermore,
$$d''Z^{q_{ic}-2\za i,p_i+2\za i}=d''Z^{2+c,p_i+2\za i}=\begin{cases}0,\;\mbox{for}\;\; c=0,\\(2a+b)^{-\za i}\lp B_{i,\za i}+
C_{i,\za i}-2D_{i,\za i}\rp z^{2+i(2\za +\zb)}\p_{12},\;\mbox{for}\;\; c=1,\\
(2a+b)^{-\za i}\lp E_{i,\za i}-F_{i,\za i}\rp z^{3+i(2\za
+\zb)}\p_{123},\;\mbox{for}\;\; c=2,\end{cases}$$ is a
$d'$-coboundary if and only if the coefficient vanishes.
\end{lem}

{\it Proof}. We must compute the differentials
$d'=\p_{\zL_{4,I}}=[\zL_{4,I},.]$ and
$d''=\p_{\zL_{4,II}}=[D^{-1}z^3Y_{12},.]=:[f\mathbf{X},.]$ on the
$Z^{q_{ic}-2k,p_i+2k}$. These cochains have the form $g{\cal
Y}:=D^{-1}X^J{\cal Y}:=D^{-1}D'^nz^m\sum_jr_j\mathbf{Y}_j,$
$n,m\in\N,r_j\in\R$, where the degree $c$ of wedge product
$\mathbf{Y}_j$ is independent of $j$. Hence, Equation
(\ref{bracketFuncWedge}) gives
$$d''(g{\cal Y})=[f\mathbf{X},g{\cal Y}]=f[\mathbf{X},g]\w{\cal Y}+(-1)^{c}g[{\cal
Y},f]\w\mathbf{X}.$$ On the other hand, Equations
(\ref{cobAdmiss}) and (\ref{Lambda1cob}) entail $d'(g{\cal
Y})=[\zL_{4,I},g{\cal Y}]=[\zL_{4,I},g]\w{\cal
Y}=\sum_{\E}X_{\E}(g)\;Y_{\E}\w{\cal Y},$ where $X_1=bY_2-aY_3,
X_2=aY_3-bY_1,X_3=a(Y_1-Y_2).$ Since
$$Y_{\E}\lp\frac{X^J}{D}\rp=(j_{\E}-1)\frac{X^J}{D}$$ (same notations as above), we get $$d'(g{\cal Y})=g\lp b(n-1)-a(m-1)\rp\lp Y_1-Y_2\rp\w{\cal
Y}.$$ In particular, we recover the result
$d'Z^{q_{ic},p_i}=ig(b\za-a\zb)(Y_1-Y_2)\w{\cal Y}=0,$ and, when
setting $a=0,b=1,{\cal Y}=1,$ we find
$$[\mathbf{X},g]=g(n-1)(Y_1-Y_2).$$\\

We now compute $d''Z^{q_{ic}-2k,p_i+2k},$
$k\in\{0,1,\ldots,\za\,i\},$ and $d'Z^{q_{ic}-2(k+1),p_i+2(k+1)},$
$k\in\{0,1,\ldots,\za\,i-1\}$.
\begin{enumerate}\item $c=0$\\

It follows from the preceding equations that
$$d''Z^{q_{i0}-2k,p_i+2k}=\complement^k_{\za i}(\za\,i-k)A_{ik}(2a+b)^{-k}D^{-1}D'^{\za i-k}z^{3+\zb i+2k}(Y_1-Y_2)$$
and that
$$d'Z^{q_{i0}-2(k+1),p_i+2(k+1)}=-\complement^{k+1}_{\za i}(k+1)A_{i,k+1}(2a+b)^{-k}D^{-1}D'^{\za
i-k}z^{3+\zb i+2k}(Y_1-Y_2).$$ Since for any $p,n\in\N,p<n$, we
have $\complement^p_n (n-p)=\complement^{p+1}_n (p+1),$ the sum of
these coboundaries vanishes if and only if $A_{i,k+1}=A_{ik},$ for
any $k\in\{0,1,\ldots,\za\,i-1\}.$ For $k=\za\,i$, we get
$$d''Z^{q_{i0}-2\za i,p_i+2\za i}=d''Z^{2,p_i+2\za i}=0.$$ \item
$c=1$\\

A short computation shows that
$$\begin{array}{rl}d''Z^{q_{i1}-2k,p_i+2k}=&\complement^{k}_{\za i}(2a+b)^{-k}D^{-1}D'^{\za i-k}z^{3+\zb
i+2k}[
-D_{ik}(\za\,i-k)Y_{23}-D_{ik}(\za\,i-k)Y_{31}\\
&+\lp(B_{ik}+C_{ik})(\za\,i-k+1)-2D_{ik}\rp Y_{12}]\end{array}$$
and that
$$\begin{array}{rl}d'Z^{q_{i1}-2(k+1),p_i+2(k+1)}=&-\complement^{k+1}_{\za
i}(k+1)(2a+b)^{-k}D^{-1}D'^{\za\,i-k}z^{3+\zb\,i+2k}[-D_{i,k+1}Y_{23}\\&-D_{i,k+1}Y_{31}+
(B_{i,k+1}+C_{i,k+1})Y_{12}].\end{array}$$ If
$k\in\{0,1,\ldots,\za\,i-1\}$, the sum of these coboundaries
vanishes if and only if $$B_{i,k+1}+C_{i,k+1}=\frac{\lp\za
i-k+1\rp \lp B_{ik}+C_{ik}\rp-2D_{ik}}{\za
i-k}\quad\mbox{and}\quad D_{i,k+1}=D_{ik}.$$ Furthermore, for
$k=\za\,i$, the first of the preceding ``coboundary equations''
provides the announced result for $d''Z^{q_{i1}-2\za i,p_i+2\za
i}$. As $\R[[z]]\p_{12}$ is part of the cohomology of
$\overline{d}_0=d'$, this $d''$-coboundary is a $d'$-coboundary if
and only if its coefficient vanishes. \item
$c=2$\\

We immediately obtain
$$d''Z^{q_{i2}-2k,p_i+2k}=\complement^{k}_{\za
i}(\za\,i-k+1)(2a+b)^{-k}D^{-1}D'^{\za i-k}z^{3+\zb
i+2k}(E_{ik}-F_{ik})Y_{123}$$ and
$$d'Z^{q_{i2}-2(k+1),p_i+2(k+1)}=-\complement^{k+1}_{\za
i}(k+1)(2a+b)^{-k}D^{-1}D'^{\za\,i-k}z^{3+\zb\,i+2k}(E_{i,k+1}-F_{i,k+1})Y_{123}.$$
Hence the announced upshots.  \rule{1.5mm}{2.5mm}
\end{enumerate}

Let us recall that the admissible values of $i$ (and the potential
conventions on coefficients $A_{10},B_{10},C_{10}$) depend on
quotient $b/a$. Moreover, for $b\slash
a\in\R^*\setminus\Q^*_+,(b,a)\nsim (-1,\za),\za\in\N^*,$ we set
$\za=1\in\N^*.$ Actually, in this case, $\za$ needed not be
defined before, as it was systematically multiplied by $i=0.$\\

The following theorem provides the complete description of the
considered SpecSeq.

\begin{theo}\label{basicTheo} The even terms $E_{2(n-1)\za+4}=E_{2(n-1)\za+6}=\ldots =E_{2n\za+2}$ $($$n\in\N$; for $n=0$,
this package contains only term $E_2$$)$ of the above defined
SpecSeq are canonically isomorphic $($i.e.
$d_{2(n-1)\za+4}=d_{2(n-1)\za+6}=\ldots =d_{2n\za}=0$$)$ and admit
the below system of representatives: \begin{enumerate}\item All
representatives of type $2$ of $E_2\sim\, ^0\!H({\cal P})$, except
$$\R z^{i(2\za+\zb)+2}\p_{12}\quad\mbox{and}\quad\R z^{i(2\za+\zb)+3}\p_{123},$$  for all
admissible $i\in\{0,1,\ldots,n-1\}$. \item All representatives of
type $1$ of $E_2\sim\, ^0\!H({\cal P})$, altered as follows:
\begin{itemize}\item For all admissible $i\in\{n,n+1,\ldots\}$, $$Z^{q_{ic},p_i}\rightsquigarrow \sum_{k=0}^{\za\, n}Z^{q_{ic}-2k,p_i+2k},$$
where the coefficients $A_{ik},B_{ik},C_{ik},D_{ik},E_{ik},F_{ik}$
incorporated into the terms of the {\small RHS} verify conditions
$(C_0)-(C_2)$ of Lemma \ref{Lemma} up to $k=\za\,n-1$. \item For
all admissible $i\in\{0,1,\ldots,n-1\}$,
$$Z^{q_{ic},p_i}\rightsquigarrow \begin{cases}\lp D'+\frac{z^2}{2a+b}\rp^{\za i}z^{\zb i}A_{i0},\quad\mbox{if}\quad c=0,\\
\lp D'+\frac{z^2}{2a+b}\rp^{\za i}z^{\zb i}\lp
B_{i0}(Y_1+\frac{1}{2}Y_3)+C_{i0}(Y_2+\frac{1}{2}Y_3)\rp,\quad\mbox{if}\quad
c=1,\\ D'^{\za i}z^{\zb
i}E_{i0}(Y_{23}+Y_{31}),\quad\mbox{if}\quad
c=2.\end{cases}$$\end{itemize}\end{enumerate}\end{theo}

{\it Proof}. The proof is by induction on $n$. For $n=0$, Theorem
\ref{basicTheo} is obviously valid. Assume now that it holds true
for $0,1,\ldots,n-1$ ($n\in\N^*$). We first transfer the
description of $E_{2(n-2)\za+4}=\ldots=E_{2(n-1)\za+2}$ to the
concrete model side, in order to compute
$\overline{d}_{2(n-1)\za+2}$. When having a look at the packages
of terms that are known to be isomorphic, we see that the only
differentials (under $\overline{d}_{2(n-1)\za+2}$) that do not
vanish are $\overline{d}_{2m\za+2}$ ($m\in\{0,1,\ldots,n-2\}$).
Hence the target of vector space isomorphism
$$I_{2(n-1)\za+2}:E_{2(n-1)\za+2}\raa\,
^{(n-2)\za+1}\!H(^{(n-3)\za+1}\!H(\ldots ^1\!H(^0\!H({\cal
P})))),$$ which---as it appears from its general
description---maps the system of $E_{2(n-1)\za+2}$-representatives
onto the system evidently made up by: \begin{enumerate}\item All
representatives of type $2$ of $E_2$, except $\R
z^{i(2\za+\zb)+2}\p_{12}$ and $\R z^{i(2\za+\zb)+3}\p_{123},$ for
all admissible $i\in\{0,1,\ldots,n-2\}$.\item All representatives
of type $1$ of $E_2$, $Z^{q_{ic},p_i},$ $i$ admissible,
$c\in\{0,1,2\}$, with, for all admissible
$i\in\{0,1,\ldots,n-2\}$, $B_{i0}+C_{i0}=2D_{i0},$ if $c=1$, and
$E_{i0}=F_{i0},$ if $c=2$.\end{enumerate}

We now compute the cohomology of space $(^{(n-2)\za+1}\!H(\ldots
^0\!H({\cal P})),\overline{d}_{2(n-1)\za+2})$. If $z^{qp}$ is one
of the representatives of the preceding system,
\be\overline{d}_{2(n-1)\za+2}[[z^{qp}]_0\ldots]_{(n-2)\za+1}=[[d''z^{q-2\za(n-1),p+2\za(n-1)}]_0\ldots]_{(n-2)\za+1},\label{cobMainTheo}\ee
where $z^{q-2\za(n-1),p+2\za(n-1)}$ is the last entry of an
arbitrary solution of $S(z^{qp};\za (n-1)+1)$.

The $d''$-coboundary of any $z^{qp}$ of type $2$ vanishes. This is
obvious if $z^{qp}$ is a $3$-cochain or has the form
$\op{Cas}(\zL_{4,I})Y_{12}$ (as
$d''=[\zL_{4,II},.]=[D^{-1}z^3Y_{12},.]$). If $z^{qp}$ is a
$2$-cochain with singularities, e.g. $D^{-1}p(x)Y_{23}$, where
$p(x)$ is a polynomial in $x$, we get
$d''z^{qp}=[D^{-1}z^3Y_{12},D^{-1}p(x)Y_{23}]=-z^3p(x)D^{-2}Y_{123}+z^3p(x)D^{-2}Y_{123}=0.$
Hence, for any type $2$ representative $z^{qp}$, system
$S(z^{qp};s)$ admits solution $(z^{qp},0,\ldots,0)$, for any
$s\in\N^*$ ($\mathbf{S_1}$, representative extended by 0
[reference needed in the following]), and Coboundary
(\ref{cobMainTheo}) vanishes.

Let now $z^{qp}$ be a representative $Z^{q_{ic},p_i}$ of the first
type. We know from Lemma \ref{Lemma} that $Z^{q_{ic}-2k,p_i+2k}$,
$k\in\{0,1,\ldots,\za\,i\}$, with coefficients that verify
($C_0$)-($C_2$), is a solution of $S(Z^{q_{ic},p_i};\za\,i+1)$.
\begin{enumerate}\item For any admissible $i\in\{n,n+1,\ldots\}$,
this solution can be truncated to a solution of
$S(Z^{q_{ic},p_i};\alpha\,n+1)$ ($\mathbf{S_2}$, truncated
standard solution). Hence, Coboundary (\ref{cobMainTheo})
vanishes. \item If $i$ is admissible in $\{0,1,\ldots,n-2\}$, we
have $$B_{i0}+C_{i0}=2D_{i0}\quad\mbox{and}\quad E_{i0}=F_{i0}.$$
It then follows from ($C_1$) and ($C_2$) that the same relation
holds for $k={\za\,i}$, i.e. that
$B_{i,\za\,i}+C_{i,\za\,i}=2D_{i,\za\,i}$ and
$E_{i,\za\,i}=F_{i,\za\,i}$. This however implies that
$d''Z^{q_{ic}-2\za\,i,p_i+2\za\,i}=0$, so that system
$S(Z^{q_{ic},p_i};\za\,n+1)$ admits an obvious solution
($\mathbf{S_3}$, standard solution extended by 0) and that
Coboundary (\ref{cobMainTheo}) vanishes again. \item If $i=n-1$ is
admissible,
$$\overline{d}_{2(n-1)\za+2}[[Z^{q_{n-1,c},p_{n-1}}]_0\ldots]_{(n-2)\za+1}=[[d''Z^{q_{n-1,c}-2\za(n-1),p_{n-1}+2\za(n-1)}]_0\ldots]_{(n-2)\za+1}.$$
In view of Lemma \ref{Lemma}, this class vanishes for $c=0$, and
coincides, if $c=1$ (resp. $c=2$), up to a coefficient, with class
$[[z^{(n-1)(2\za +\zb)+2}\p_{12}]_0\ldots]_{(n-2)\za+1}$ (resp.
$[[z^{(n-1)(2\za +\zb)+3}\p_{123}]_0\ldots]_{(n-2)\za+1}$). The
above depicted system of representatives of
$^{(n-2)\za+1}\!H(\ldots ^0\!H({\cal P}))$ shows that the
preceding two classes do not vanish. Hence, the cocycle-condition
is equivalent with the annihilation of the mentioned coefficient,
i.e. with
$$B_{n-1,\za (n-1)}+C_{n-1,\za (n-1)}=2D_{n-1,\za (n-1)}\quad
(\mbox{resp.}\quad E_{n-1,\za (n-1)}=F_{n-1,\za (n-1)}),$$ or, as
already explained, \be B_{n-1,0}+C_{n-1,0}=2D_{n-1,0}\quad
(\mbox{resp.}\quad E_{n-1,0}=F_{n-1,0}).\label{NewCoeffCond}\ee
Since it clearly follows from our computations that the space of
$\overline{d}_{2(n-1)\za+2}$-coboundaries is generated by the two
just encountered non-vanishing classes, cohomology space
$^{(n-1)\za+1}\!H(^{(n-2)\za+1}\!H(\linebreak\ldots ^0\!H({\cal
P})))$ has the same system of representatives than its predecessor
$^{(n-2)\za+1}\!H(\ldots ^0\!H({\cal P}))$, but with exclusions
carried out and conditions on $B_{i0},C_{i0},D_{i0},E_{i0},F_{i0}$
valid for all admissible $i\in\{0,1,\ldots,n-1\}.$
\end{enumerate}

It now suffices to apply Corollary \ref{Corollary} to cohomology
space $^{(n-1)\za+1}\!H(^{(n-2)\za+1}\!H(\ldots ^0\!H({\cal
P}))).$ Observe first that ($\mathbf{S_1}$)-($\mathbf{S_3}$)
entail existence of a solution of $S(z^{qp};\za\,n+1)$, for all
representatives $z^{qp}$ dissimilar from $Z^{q_{n-1,c},p_{n-1}}$.
But, as the coefficients of these last representatives---viewed as
representatives of the preceding
$\overline{d}_{2(n-1)\za+2}$-cohomology space---satisfy Conditions
(\ref{NewCoeffCond}), the coboundaries
$d''Z^{q_{n-1,c}-2\za(n-1),p_{n-1}+2\za(n-1)}$ vanish. So the
previously met solution of $S(Z^{q_{n-1,c},p_{n-1}};\za\,(n-1)+1)$
can be indefinitely extended by $0$ ($\mathbf{S_4}$, standard
solution extended by 0). Finally, Corollary \ref{Corollary} is
applicable for $s=\za\,n+1$.

Hence, spaces $E_{2(n-1)\za+4}=\ldots=E_{2n\za+2}$ coincide and we
build, from the known system $z^{qp}$ of representatives of
$^{(n-1)\za+1}\!H(\ldots ^0\!H({\cal P}))$, a system of
$E_{2n\za+2}$ by just summing-up the entries of any solutions of
the systems $S(z^{qp};\za\,n+1)$. For any $Z^{q_{ic},p_i}$, the
coefficients of which verify
$$B_{i0}+C_{i0}=2D_{i0}\quad\mbox{(}c=1\mbox{)}\quad\mbox{and}\quad
E_{i0}=F_{i0}\quad\mbox{(}c=2\mbox{)},$$ the standard
$Z^{q_{ic}-2k,p_i+2k}$, $k\in\{0,1,\ldots,\za\,i\},$ are solution,
see Lemma \ref{Lemma}, of $S(Z^{q_{ic},p_i};\za\,i+1)$, e.g. if we
choose \be A_{ik}=A_{i0}\;\,\mbox{(}c=0\mbox{)},\quad
B_{ik}=B_{i0},C_{ik}=C_{i0},D_{ik}=D_{i0}\;\,\mbox{(}c=1\mbox{)},\quad\mbox{and}\quad
E_{ik}=F_{ik}=0\;\,\mbox{(}c=2,k\neq
0\mbox{)}.\label{CoeffChoice}\ee If we pull the concrete side
representatives back to theoretical side representatives using
these solutions, we exactly get, see $S_1$-$S_4$, the sought-after
system.
\rule{1.5mm}{2.5mm}\\

{\bf Remark}. \begin{enumerate} \item We already observed
previously the obvious fact that when pulling {\small
RHS}-representatives back, using different solutions of the
standard system, we obtain equivalent {\small
LHS}-representatives. These equivalent {\small
LHS}-representatives would implement cohomologous cocycles in
cohomology space $H(\zL_4)$. Choice (\ref{CoeffChoice}) will
induce in cohomology the most basic possible cocycles. \item Note
also that in view of Theorem \ref{basicTheo} and our conventions
on coefficients $B_{10},C_{10}$, cocycle $\R{\cal A}_{\za}Y_3$
disappears from all spaces $E_{2r}$, $r\ge 2\za+4$.\end{enumerate}

\subsection{Limit of the SpecSeq and reconstruction of the cohomology}

The limit of the SpecSeq can be guessed from Theorem
\ref{basicTheo}. However, we already stressed the importance of a
careful reading of all results through the isomorphisms involved
in the theory of spectral sequences. The proof of Theorem
\ref{basicTheo} shows for instance that the appropriate Casimir
functions appear, when we pull the {\small RHS}-representatives
back to the {\small LHS}, i.e. read them through isomorphism
$I^{-1}_{2(n-1)\za+3}$. Hence, a precise description of the
isomorphisms that lead now to the cohomology of $\zL_4$ is
essential.

\subsubsection{General results}

Let us consider the SpecSeq associated with a (regular) filtered
graded differential space ($K,d,K_p,\linebreak K^n$) and recall
that the limit spaces
$E_{\infty}^{pq},Z_{\infty}^{pq},B_{\infty}^{pq}$ are defined
exactly as spaces $E_r^{pq},Z_r^{pq},B_r^{pq}$, see Section
\ref{ModSpecSeq}, so that $Z_{\infty}^{pq}$ and $B_{\infty}^{pq}$
are the spaces of cocycles and coboundaries in $K_p\cap K^{p+q}$
respectively. For any fixed $p$ and $q$, regularity implies that
the target space of the restriction of $d_r$ to $E_r^{pq}$
vanishes, if $r>q+1$. Thus, there is a canonical linear surjective
map $\vartheta_r^{pq}:E_r^{pq}\raa H^{pq}(E_r,d_r)\raa
E_{r+1}^{pq}$. For $s\ge r>q+1$, we define
$\zy_{rs}^{pq}:=\vartheta_{s-1}^{pq}\circ\ldots\circ\vartheta_r^{pq}:E_r^{pq}\raa
E_s^{pq},$ and for $r>q+1$, we set \be\zy_r^{pq}:E_r^{pq}\ni
[\frak{z}_r^{pq}]_{E_r^{pq}}\raa
[\frak{z}_r^{pq}]_{E_{\infty}^{pq}}\in
E_{\infty}^{pq}.\label{IsomInducLim}\ee Due to regularity, the
first two of the well-known inclusions $Z_{\infty}^{pq}\subset
Z_r^{pq}, Z_{\infty}^{p+1,q-1}\subset Z_{r-1}^{p+1,q-1},$ and
$B_{r-1}^{pq}\subset B_{\infty}^{pq}$ are actually double
inclusions, and $Z_{\infty}^{p+1,q-1}+B_{r-1}^{pq}\subset
Z_{\infty}^{p+1,q-1}+B_{\infty}^{pq}\subset Z_{\infty}^{pq}$.
Hence, map $\zy_r^{pq}$ is canonical, linear and surjective. It is
known that space $E_{\infty}^{pq}$ together with the preceding
linear surjections $\zy_r^{pq}$ is a model of the inductive limit
of the inductive system $(E_r^{pq},\zy_{rs}^{pq}).$ Consider now a
first quadrant SpecSeq (i.e. $p,q\in\N$) and assume that $K_0=K$.
For any $p,q$, the SpecSeq collapses at $$r>\op{sup}(p,q+1),$$
more precisely, $E_r^{pq}=E_{\infty}^{pq}$ and
$\zy_r^{pq}=\op{id}$. Indeed, in this case, in addition to the
aforementioned double inclusions ($r>q+1$), we now have also
$B_{r-1}^{pq}=K_p\cap dK_{p+1-r}\cap K^{p+q}=K_p\cap dK_0\cap
K^{p+q}=B_{\infty}^{pq}$ ($r>p$). Hence the announced results.\\

The SpecSeq associated with any filtered graded differential space
is convergent in the sense that limit $E_{\infty}^{pq}$ is known
to be isomorphic as a vector space with term $G^{pq}$ of the
bigraded space $G(H(K))$, $G$ for short, associated with the
filtered graded space $H(K)$. Let us recall that the filtration of
$H(K)$ is induced by that of $K$. More precisely, injection
$i:(K_p,d)\raa (K,d)$ is a morphism of differential spaces and
$H_p:=i_{\sharp}H(K_p)\subset H(K)$ is the mentioned filtration of
$H(K)$. In order to reduce notations, we denote the terms of the
grading of $H(K)$ simply by $H^n$. It is a fact that the
filtration and the grading of $H(K)$ are compatible and that
filtration $H_p$ is regular if its generatrix $K_p$ is. Hence,
$H_p=\oplus_{q\in\N} H_p\cap H^{p+q}=:\oplus_{q\in\N}H^{p+q}_p.$
Finally, it is a matter of knowledge that the isomorphism, say
$\zi$, between $G^{pq}:=H_p^{p+q}\slash H_{p+1}^{p+q}$ and
$E_{\infty}^{pq}$ is canonical, \be\zi: E_{\infty}^{pq}\ni
[\frak{z}_{\infty}^{pq}]_{E_{\infty}^{pq}}\raa
[[\frak{z}_{\infty}^{pq}]_{H_p^{p+q}}]_{G^{pq}}\in
G^{pq}.\label{IsomGradSpace}\ee \linebreak[1]

We now reconstruct $H(K)$ from $G$. Let us again focus on a first
quadrant SpecSeq associated with a (regular) filtered complex
($K,d,K_p,K^n$) (such that $K_0=K$). For any $n\in\N$, we denote
by $G^{n-j_1,j_1},G^{n-j_2,j_2},\ldots,G^{n-j_{k_n},j_{k_n}},$
$n\ge j_1>j_2>\ldots >j_{k_n}\ge 0$, the non vanishing
$G^{pq}=H_p^{p+q}\slash H_{p+1}^{p+q}$, $p+q=n.$ Since $H_0=H(K)$
and $H_{p}^n=H_p\cap H^n=0,\forall p>n$, it follows that
\begin{eqnarray*}H^n=H_0^n=\ldots =H_{n-j_1}^n\supset H_{n-j_1+1}^n=\ldots
=H_{n-j_2}^n\\\supset H_{n-j_2+1}^n\ldots H_{n-j_{k_n}}^n\supset
H_{n-j_{k_n}+1}^n= \ldots =H_n^n=0.\end{eqnarray*} Hence,
$$ H^n\slash H_{n-j_2}^n=G^{n-j_1,j_1},\ldots, H_{n-j_{k_n-1}}^n\slash
H_{n-j_{k_n}}^n=G^{n-j_{k_n-1},j_{k_n-1}},H_{n-j_{k_n}}^n=G^{n-j_{k_n},j_{k_n}}.$$

However, if $B\slash A=C$, $A$ a vector subspace of $B$, the
sequence $0\raa A\stackrel{i}{\raa}B\stackrel{p}{\raa}C\raa 0,$ is
a short exact sequence of vector spaces. A short exact sequence in
a category is split if and only if kernel $A$ admits in vector
space $B$ a complementary subspace that is a subobject,
or---alternatively---if and only if there is a right inverse
morphism $\chi:C\raa B$ of projection $p$. Of course, in the
category of vector spaces such a sequence is always split. If
$\chi$ is a linear right inverse of $p$, we have
$B=A\oplus\chi(C)$.

Let us now come back to our circumstances. If
$\chi_{1},\ldots,\chi_{k_n-1}$ denote splittings of the involved
sequences, central extension $H^n$ is given by \be
H^n=\chi_{1}(G^{n-j_1,j_1})\oplus\ldots\oplus\chi_{k_n-1}(G^{n-j_{k_n-1},j_{k_n-1}})\oplus
G^{n-j_{k_n},j_{k_n}}.\label{ReconstCohom}\ee It follows of course
from Equation (\ref{ReconstCohom}) that $H(K)$ is---in this vector
space setting---isomorphic with $G=G(H(K)).$ It is known that in
the case of ring coefficients, extension problems may prevent the
reconstruction of $H(K)$ from $G(H(K)).$

\subsubsection{Application to Poisson tensor $\mathbf{\zL_4}$}

The next proposition provides a system of representatives of the
cohomology space of
$$\zL_4=ayz\p_{23}+axz\p_{31}+(bxy+z^2)\p_{12}\quad (a\neq 0,b\neq 0).$$ Remember that
$D'=xy$ and $Y_1=x\p_1,Y_2=y\p_2,Y_3=z\p_3$. If
$\frac{b}{a}\sim\frac{\zb}{\za}\in\Q^*_+,$ we define
$$\op{Cas}(\zL_4):=\oplus_{i\in\N}\R\lp
D'+\frac{z^2}{2a+b}\rp^{\za\,i}z^{\zb\,i}$$ and use the above
introduced notation $\op{Cas}(\zL_{4,I})=\oplus_{i\in\N}\R
D'^{\za\,i}z^{\zb\,i}.$ If $\frac{b}{a}\in\R^*\setminus\Q^*_+,$ we
set $\op{Cas}(\zL_4):=\R$ and, as aforementioned, ${\cal
A}_{\za}=D'^{\za}z^{-1}$.

\begin{theo}

\begin{enumerate}
\item If $\frac{b}{a}\in\Q^*_+,$ the cohomology of $\zL_4$ is
given by
\begin{eqnarray*} E_{\infty}\sim G\sim H(\zL_4)&=&\op{Cas}(\Lambda
_{4})\oplus\op{Cas}(\Lambda
_{4})(Y_{1}+\frac{1}{2}Y_3)\oplus\op{Cas}(\Lambda
_{4})(Y_2+\frac{1}{2}Y_3)\\&&\oplus \op{Cas}(\Lambda
_{4,I})(Y_{23}+Y_{31})\oplus\op{Cas}(\zL_{4,I})Y_{12}
\oplus \op{Cas}(\Lambda _{4,I})Y_{123}\\
&&\oplus\bigoplus_{k\in\N\setminus\N\,(2\za+\zb)+2}\R z^k\partial
_{12}\oplus\bigoplus_{k\in\N\setminus\N\,(2\za+\zb)+3}\R
z^k\partial
_{123}\\
&&\oplus \left\{
\begin{array}{l}
\R[[x]]\p_{23}\oplus\R[[y]]\p_{31}\oplus
(\R[[x]]\oplus\R[[y]])\partial _{123},\text{if }b=a\\
0,\text{otherwise}
\end{array}
\right.
\end{eqnarray*}

\item If $\frac{b}{a}\in\R^*\setminus\Q^*_+,$ we have
\begin{eqnarray*}
E_{\infty}\sim G\sim H(\zL_4)&=&\op{Cas}(\Lambda
_{4})\oplus\op{Cas}(\Lambda
_{4})(Y_1+\frac{1}{2}Y_3)\oplus\op{Cas}(\zL_4)(Y_2+\frac{1}{2}Y_3)\\&&\oplus
\op{Cas}(\Lambda
_{4})(Y_{23}+Y_{31})\oplus\op{Cas}(\zL_4)Y_{12}\oplus
\op{Cas}(\Lambda _{4})Y_{123}\\&&\oplus \left\{
\begin{array}{l}\oplus\R{\cal A}_{\za}(Y_{23}+Y_{31})\oplus \R{\cal
A}_{\za}Y_{123},\text{ if }(b,a)\sim(-1,\za)\\0,\text{otherwise}
\end{array}
\right.\\&&\oplus
\begin{cases}\bigoplus_{k\in\N\setminus\{2,2\za+1\}}\R z^k\partial
_{12}\oplus\bigoplus_{k\in\N\setminus\{3,2\za+2\}}\R z^k\partial
_{123},\text{ if
}(b,a)\sim(-1,\za)\\\bigoplus_{k\in\N\setminus\{2\}}\R z^k\partial
_{12}\oplus\bigoplus_{k\in\N\setminus\{3\}}\R z^k\partial
_{123},\text{otherwise}\end{cases}
\end{eqnarray*}
\end{enumerate}
\label{IndLim}\end{theo}

{\it Proof}. Fix $a,b\in\R^*$ and take any representative of
$E_2$. Remember that the representatives of type 1 are exactly the
cochains $Z^{q_{ic},p_i}$ ($i$ admissible, $c\in\{0,1,2\}$).
Moreover, we say that a representative of type 2 is critical if it
has the form $\R z^{i(2\za+\zb)+2}\p_{12}$ or $\R
z^{i(2\za+\zb)+3}\p_{123}$ ($i$ admissible). If the considered
representative $z^{qp}$ is of type 2 and not critical (resp. of
type 2 and critical, of type 1), we choose $n\in\N$ such that
$2n\za+2>\op{sup}(p,q+1)$ (resp. $2n\za+2>\op{sup}(p,q+1,2i\za+2)$
[hence, we have $n-1\ge i$],
$2n\za+2>\op{sup}(p_i,q_{ic}+1,2i\za+2)$). The system of
representatives of $E_{\infty}\sim G$ specified in Theorem
\ref{IndLim} arises now from Theorem \ref{basicTheo} and from the
canonical isomorphisms (\ref{IsomInducLim}) (condition:
$r>\op{sup}(p,q+1)$) and (\ref{IsomGradSpace}). These
representatives are representatives of bases of the non vanishing
$G^{pq}$. In order to compute $H(K),$ it suffices to build
arbitrary splittings in keeping with Equation
(\ref{ReconstCohom}). Hence, it suffices to choose, for any class
of any basis of the concerned $G^{pq},$ an arbitrary
representative, e.g. the aforementioned one. It follows that
$H(K)$ admits exactly the same representatives as $E_{\infty}\sim G$. \rule{1.5mm}{2.5mm}\\

{\bf Remarks}. Hence, the twist makes a threefold impact on
cohomology. When applying our computing device, see Theorem
\ref{basicTheo}, we get, at {\it each} turn of the handle, on the
model level, roughly speaking, cocycle-conditions on the
coefficients related with an {\it additional} power of the basic
Casimir $C_{\zL_{4,I}}$ of $\zL_{4,I}$, and we exclude a {\it
supplementary} pair of singularity-induced classes. These
conditions appear in cohomology as terms $Y_i$ or $Y_{ij}$ with
the same ``Casimir-coefficient''. Eventually, the
cocycle-conditions allow to lift the mentioned {\it accessory}
power of $C_{\zL_{4,I}}$ to the real level as power of Casimir
$C_{\zL_4}$ of $\zL_4$ or---depending on cochain degree---as power
of Casimir $C_{\zL_{4,I}}$. We know that such a lift is not unique
and that two different ones are cohomologous. It follows from
Theorem \ref{IndLim} (resp. from the proof of Theorem
\ref{basicTheo}) that any term of $\op{Cas}(\zL_4)(Y_{23}+Y_{31})$
is a $\zL_4$-cocycle (resp. can be chosen as lift of the
corresponding term in $\op{Cas}(\zL_{4,I})(Y_{23}+Y_{31})$, as
well as this term itself). So any term of
$\op{Cas}(\zL_4)(Y_{23}+Y_{31})$ is cohomologous to the analogous
term in $\op{Cas}(\zL_{4,I})(Y_{23}+Y_{31})$. Finally, the
aforementioned proof allows to see that $\zL_{4,I}$-cocycle
$\R{\it A}_{\za}Y_3=\R D'^{\za}z^{-1}Y_3$, which is not a product
of two $\zL_4$-cocycles, is a $\zL_4$-cocycle if and only if its
coefficient vanishes.\\

Let us in the end have a look at singularities. The singular locus
of $\zL_{4,I}$ (resp. $\zL_4$) is made up by the three coordinate
axes (resp. the axis of absciss\ae\, and the axis of ordinates).
Comparing the results of Proposition \ref{cohoZL4I} and of Theorem
\ref{IndLim}, we see that the twist $\zL_{4,II}$, which removes
the $z$-axis from the singular locus, cancels only part of the
corresponding polynomials in cohomology. We already observed in
\cite{MP} that, for $r$-matrix induced tensors, some coefficients
of non bounding $2$- or $3$-cocycles can just be interpreted as
polynomials on singularities via an extension of the polynomial
ring of the singular locus. In the case of twisted $r$-matrix
induced structures, some of these polynomial coefficients are
simply not polynomials on singularities.

\section{Formal cohomology of Poisson tensor $\zL_8$}

We now describe the cohomology space of the twisted quadratic
Poisson structure $$\zL_8=\lp\frac{{\frak a}+{\frak
b}}{2}(x^2+y^2)\pm z^2\rp\p_{12}+{\frak a}xz\p_{23}+{\frak
a}yz\p_{31}\quad ({\frak a}\neq 0, {\frak b}\neq 0).$$ If we
substitute $c$ (resp. $b$) for $-{\frak b}$ (resp. $({\frak
a}+{\frak b})/2$), tensor $\zL_8$ reads
\be\zL_{8}=b(x^2+y^2)\p_{12}+(2b+c)xz\p_{23}+(2b+c)yz\p_{31}\pm
z^2\p_{12}.\label{lambda8}\ee Henceforth we use parameters $b$ and
$c$. Assumptions ${\frak a}\neq 0, {\frak b}\neq 0$ are equivalent
with $2b+c\neq 0, c\neq 0$. Moreover, the $r$-matrix induced part
$\zL_{8,I}=b(x^2+y^2)\p_{12}+(2b+c)xz\p_{23}+(2b+c)yz\p_{31}$ of
$\zL_8$ is nothing but structure $\zL_7$ with parameter $a=0$, see
\cite[Section 9]{MP}, so that term $E_2\simeq H(\zL_{8,I})$ of the
spectral sequence follows from \cite[Theorems 6,8,9]{MP}.\\

Let us recall that the $Y_i$ stem from $\zL_{8,I}$, i.e. from
$\zL_7$. Hence, $Y_1=x\p_1+y\p_2,Y_2=x\p_2-y\p_1,Y_3=z\p_3$. We
set $D'=x^2+y^2$. Moreover, if $\frac{b}{c}\in\Q, b(2b+c)<0,$ we
denote by $(\zb,\zg)\sim (b,c)$ the irreducible representative of
the rational number $\frac{b}{c}$, with positive denominator,
$\zb\in\Z,\zg\in\N^*,$ and if $\frac{b}{c}\in\Q, b(2b+c)>0,$
$(\zb,\zg)\sim (b,c)$ denotes the irreducible representative with
positive numerator, $\zb\in\N^*,\zg\in\Z^*.$

\begin{theo}\label{IndLim8}
The terms of the cohomology space of $\Lambda_{8}$ (see
(\ref{lambda8})) are given by the following equations:
\begin{enumerate}
\item If $\frac{b}{c}\in \mathbb{Q},b(2b+c)>0,$
\begin{eqnarray*}
H^{0}(\Lambda _{8})
&=&\op{Cas}(\Lambda_{8})=\oplus_{i\in\mathbb{N},\zg i\in
2\Z}\mathbb{R} \left(D^{\prime
}\pm\frac{z^{2}}{3b+c}\right) ^{(\zb +\frac{\zg}{2})i}z^{\zb\,i}\,,\\
H^{1}(\Lambda _{8}) &=&\op{Cas}(\Lambda _{8,I})Y_{2}\oplus
\op{Cas}(\Lambda_{8})(Y_{1}+Y_{3})\,,\\
H^{2}(\Lambda _{8}) &=&\op{Cas}(\Lambda _{8,I})Y_{12}\oplus
\op{Cas}(\Lambda _{8,I})Y_{23}\oplus \bigoplus_{\substack{
k\in\N\setminus\N(3\zb +\zg)+2}}\mathbb{R}z^{k}\partial _{12}\,,\\
H^{3}(\Lambda _{8}) &=&\op{Cas}(\Lambda _{8,I})Y_{123}\oplus
\bigoplus_{\substack{ k\in\N\setminus\N(3\zb +\zg )+3}} \mathbb{R}
z^{k}\partial _{123}\,,\end{eqnarray*} where $\op{Cas}(\Lambda
_{8,I})=\oplus _{i\in\mathbb{N},\zg i\in 2\Z} \mathbb{R} D^{\prime
(\zb +\frac{\zg }{2})i}z^{\zb i}.$

\item If $\frac{b}{c}\notin\Q$ or $\frac{b}{c}\in
\mathbb{Q},b(2b+c)<0$,
\begin{eqnarray*}
H^{0}(\Lambda _{8}) &=&\op{Cas}(\Lambda _{8})=\mathbb{R}\,,\\
H^{1}(\Lambda _{8}) &=&\op{Cas}(\Lambda _{8})Y_{2}\oplus
\op{Cas}(\Lambda
_{8})(Y_{1}+Y_{3})\,,\\
H^{2}(\Lambda _{8}) &=&\op{Cas}(\Lambda _{8})Y_{12}\oplus
\op{Cas}(\Lambda _{8})Y_{23}\\&&\oplus
\begin{cases}\bigoplus_{k\in\mathbb{N}\setminus \{2,\zg-1\}}\mathbb{R}
z^{k}\partial _{12},\text{if }(b,c)\sim (-1,\zg),\;\zg \in
\left\{4,6,8,...\right\}\\ \bigoplus_{k\in\mathbb{N}\setminus
\{2\}}\mathbb{R}
z^{k}\partial _{12},\text{otherwise }\end{cases} \\
&&\oplus \begin{cases} \mathbb{R} {\cal A}_{\zg }Y_{23},\text{if
}(b,c)\sim (-1,\zg ),\zg \in \left\{4,6,8,...\right\}\\
0,\text{otherwise ,}
\end{cases}\\
H^{3}(\Lambda _{8}) &=&\op{Cas}(\Lambda _{8})Y_{123}\\&&\oplus
\begin{cases}\bigoplus_{k\in\N\setminus\{3,\zg\}}\mathbb{R} z^{k}\partial
_{123},\text{if }(b,c)\sim (-1,\zg ),\zg\in \left\{
4,6,8,...\right\}\\\bigoplus_{k\in\N\setminus\{3\}}\mathbb{R}
z^{k}\partial
_{123},\text{otherwise}\end{cases}\\
&&\oplus \begin{cases} \mathbb{R} {\cal A}_{\zg }Y_{123},\text{if
 }(b,c)\sim (-1,\zg),\zg \in
\left\{ 4,6,8,....\right\}\\ 0,\text{otherwise ,}
\end{cases}\end{eqnarray*}
where ${\cal A}_{\zg}=D^{\prime \frac{\zg}{2}-1}z^{-1}.$

\item If $b=0,$
\begin{eqnarray*} H^{0}(\Lambda _{8})
&=&\op{Cas}(\Lambda _{8})=\oplus _{i\in \mathbb{N}}\mathbb{R}
\left( D^{\prime }\pm\frac{z^{2}}{c}\right) ^{i}\,, \\
H^{1}(\Lambda _{8}) &=&\op{Cas}(\Lambda_{8,I})Y_{2}\oplus
\op{Cas}(\Lambda
_{8})(Y_{1}+Y_{3})\,, \\
H^{2}(\Lambda _{8}) &=&\op{Cas}(\Lambda _{8,I})Y_{12}\oplus
\op{Cas}(\Lambda _{8,I})Y_{23}\oplus \bigoplus_{k\in\N\setminus\{
2\N+2\}}\mathbb{R} z^{k}\partial _{12}\,, \\
H^{3}(\Lambda _{8}) &=&\op{Cas}(\Lambda _{8,I})Y_{123}\oplus
\bigoplus_{k\in\N\setminus\{2\N+3\}} \mathbb{R} z^{k}\partial
_{123}\,,
\end{eqnarray*}
where $\op{Cas}(\Lambda _{8,I})=\oplus _{i\in
\mathbb{N}}\mathbb{R} D^{\prime i}.$
\end{enumerate}
\end{theo}

\end{document}